\documentclass[12pt,reqno]{amsart}
\usepackage{amsmath}
\usepackage{amssymb}
\usepackage{amsthm}
\usepackage{mathrsfs}
\usepackage{color}
\usepackage{graphicx}
\usepackage{bbm}
\usepackage[colorlinks,linkcolor=blue,anchorcolor=red,citecolor=blue]{hyperref}
\usepackage[margin=1in]{geometry}
\usepackage{orcidlink}
\let\a=\alpha

\newcommand\bZ{{\mathbb Z}}
\newcommand\bR{{\mathbb R}}

\newcommand\bS{{\mathbb S}}

\newcommand\cG{{\mathcal G}}

\def\dis{\displaystyle}

\newtheorem{lemm}{Lemma} 
\newtheorem{theo}{Theorem} 
\newtheorem{defi}{Definition}
\newtheorem{prop}{Proposition}
\newtheorem{rema}{Remark}

\numberwithin{equation}{section}

\usepackage{color}

\usepackage{refcheck}

\makeatletter
\@namedef{subjclassname@2020}{%
	\textup{2020} Mathematics Subject Classification}
\makeatother

\begin{document}

\title{ON THE VILLANI CONJECTURE FOR MEASURE-VALUED SOLUTIONS TO THE HOMOGENEOUS BOLTZMANN EQUATION WITH DEBYE--YUKAWA POTENTIAL}

\author[Y. Tu]{Yue Tu} 
\address[YT]{School of mathematics, Shandong University, Jinan,Shandong, P. R. China} 
\email{202311876@mail.sdu.edu.cn}

\author[S.-K. Wang]{Shuaikun Wang} 
\address[SKW]{School of mathematics, Shandong University, Jinan,Shandong, P. R. China} 
\email{skwang@email.sdu.edu.cn}


\begin{abstract}
	
    We investigate the smoothing effect of measure-valued solutions to the spatially homogeneous Boltzmann equation with Debye--Yukawa potential for Maxwellian-type cross sections. The main result of this paper is a proof of the Villani conjecture in this setting: any initial datum that is not a single Dirac mass and has finite moments of order $\alpha\in(0,2]$ gives rise to a measure-valued solution which becomes infinitely smooth in the Sobolev sense for every positive time. The key ingredients are a newly adapted coercivity estimate and a reverse Povzner inequality, which provide sufficient regularization despite the weak integrability of the angular kernel. This extends the inverse power law result in \cite{morimoto2015smoothing} and removes the additional assumptions imposed in the conditional smoothing result of \cite{wang2024existence} for the Debye--Yukawa case.
\end{abstract}
\maketitle	

\section{Introduction} 
	We consider the homogeneous Boltzmann equation
	\begin{equation}  \label{1.1}
		\partial_t f(t,v)=Q(f,f)(t,v),
	\end{equation}	
	where $f(t, v)$ is the density distribution of particles with velocity $v\in \bR^3$ at time $t$. The most important part of this equation is the collision operator with the bilinear form, and it describes the change rates of the density distribution through elastic binary collisions,
	\begin{align*}
		Q(f,f)(t,v)=\int_{\bR^3} \int_{\bS^2} B(v-v_*,\sigma) [f(v')f(v_*')-f(v)f(v_*)]d\sigma dv_*,
	\end{align*}
	where $\sigma \in \bS^2$, and
	$$
	\begin{aligned}
		v'=\frac{v+v_*}{2}+\frac{|v-v_*|}{2}\sigma,	\quad v_*'=\frac{v+v_*}{2}-\frac{|v-v_*|}{2}\sigma,
	\end{aligned}
	$$
	that follows from the conservation of momentum and energy,
	\begin{align*}
		v'+v_*'=v+v_*,\quad |v'|^2+|v_*'|^2=|v|^2+|v_*|^2.
	\end{align*}
For  {monatomic} gas, the nonnegative cross section $	B(v-v_*,\sigma)$ depends only on $|v-v_*|$ and the scalar product $ \frac{v-v_*}{|v-v_*|}\cdot \sigma=\cos\theta$, where $\theta$ is the deviation angle. In most of the cases, it cannot be made explicit. For the inverse power law potentials $U(\rho)=\rho^{-(s-1)}$, $s>2$, $B(v-v_*,\sigma)$ takes the form 
\begin{equation}\label{1.2}
    \begin{aligned}
         B(v&-v_*,\sigma)=|v-v_*|^{\gamma} b(\cos\theta),\mbox{ for }   \ \gamma\geq-3 ,\
	\cos\theta=\frac{v-v_*}{|v-v_*|}\cdot \sigma,\\
    &b(\cos \theta)\theta^{2+2s}\to K( \theta\to0),\mbox{ for } \ 0<s<1,\mbox{ and }\ K>0.
    \end{aligned}
\end{equation}
	As usual, it refers to the hard potential if $\gamma>0$, and refers to the
	soft potential if $-3<\gamma<0$, and for the case $\gamma=0$ it refers to the Maxwellian molecule. For the Debye--Yukawa potentials $U(\rho)=\rho^{-1}e^{-\rho^{s}}$ with $s>0$ (the classical Debye--Yukawa potential corresponding to $s=1$), it was proved
	in \cite{villani2002review,wang2024existence} that $B(v-v_*,\sigma)$ has the singularity in the form
	\begin{equation*}
		B(v-v_*,\sigma)\approx K|v-v_*|\theta^{-2}|\log(|v-v_*|^2\theta)^{-1}|^{\frac{2}{s}-1},
	\end{equation*}
	in particular, in the case of Maxwellian molecule, $B(v-v_*,\sigma)$ has the singularity in the form
	\begin{equation}
	B(v-v_*,\sigma)=b(\cos\theta)\approx K\theta^{-2}(\log\theta^{-1})^{\frac{2}{s}-1}.
	\label{1.3}
	\end{equation}
    Let $m=2/s-1$ with $0<s<2$, in this paper, we will present the smoothing effect
for all $m>0$, which includes the true Debye--Yukawa case $m = 1$.
	\par 
	Mathematically, the range of $\theta$ can be restricted to $[0,\frac{\pi}{2}]$, by replacing $b(\cos\theta)$ by its “symmetrized” version
	$$
	[b(\cos\theta)+b(\cos(\pi-\theta))]1_{0\leq\theta\leq\frac{\pi}{2}},
	$$ 
	because of the invariance of the product
$f(v')f(v'_{\ast})$ in the collision operator $Q(f,f)$ under the change of variables $\sigma\to-\sigma$.
	\par 
 Consider the Cauchy problem of \eqref{1.1} with initial datum
	\begin{equation}
		f(0,v)=dF_{0}(v),
		\label{1.4} 
	\end{equation}
	by normalization and simple
	scaling we can assume that the initial datum $F_{0}(v)$ is a probability measure with finite moments up to the order $\alpha$, i.e.
	$F_{0}(v)\in P_{\alpha}(\bR^3)$. Especially, for any $\alpha\in (0,2]$, we denote
		\begin{align*}
			P_{\alpha}(\bR^3)=&\{F(v)\in P_{0}(\bR^3)|\int_{\bR^3}dF(v)=1,\int_{\bR^3}|v|^{\alpha}dF(v)<+\infty, \\
			& \mbox{ and if } \alpha>1,\int_{\bR^3}v_{j}dF(v)=0,j=1,2,3 \}.
		\end{align*}
	Following Jacob \cite{jacob2001pseudo} and Cannone and Karch \cite{cannone2010infinite}, the Fourier transform of a probability
	measure $F(v)\in P_0(\bR^3)$ is defined by
	$$
	\varphi=\mathcal{F}(F)(\xi)=\int_{\bR^3}e^{-iv\cdot\xi}dF(v),
	$$
	which is called a characteristic function of $F(v)$. The Fourier transform of space $P_\alpha(\bR^3)$ 
	has been well characterized in \cite{morimoto2015new,morimoto2017moment}.
	\par 
	Put $\mathcal{K}=\mathcal{F}(P_0(\bR^3))$. Inspired by a series of works by Toscani and his co-authors \cite{carlen1999propagation,gabetta1995metrics}, Cannone and Karch defined a subspace $\mathcal{K}^\alpha$ for $\alpha \geq 0$ as follows:
	\begin{align*}
		\mathcal {K}^{\alpha}=\{\varphi \in \mathcal {K};||\varphi-1||_{\alpha}<\infty\},
	\end{align*}
	where
	\begin{align*}
		||\varphi-1||_{\alpha}=\sup\limits_{ {\xi \in \bR^3\setminus\{0\}}}\frac{|\varphi-1|}{|\xi|^\alpha}.
	\end{align*}
	The space $\mathcal{K}^{\alpha}$  endowed with the distance
	\begin{align*}
		||\varphi-\overset{\sim}{\varphi}||_{\alpha}=\sup\limits_{ {\xi \in \bR^3\setminus\{0\}}}\frac{|\varphi-\overset{\sim}{\varphi}(\xi)|}{|\xi|^\alpha}
	\end{align*}
	is a complete metric space (see Proposition 3.10 in \cite{cannone2010infinite}) and the following embeddings (Lemma 3.12 in \cite{cannone2010infinite}) hold:
	$$
	\mathcal {K}^{\alpha}\subset\mathcal {K}^{\beta}\subset\mathcal {K}^{0}=\mathcal {K}, \mbox{ for all }\ 2\geq\alpha\geq\beta\geq0.$$
	The definition of $\mathcal{K}^{\alpha}$ is natural, because $\mathcal{F}(P_{\alpha})\subseteq \mathcal{K}^{\alpha}$ (Lemma 3.15
	in \cite{cannone2010infinite}). But the reverse of this relation does not hold, in fact, the space of $\mathcal{K}^{\alpha}$ is larger than the set
	of the Fourier transform of $P_{\alpha}$. So we introduce $\tilde{P}_{\alpha}(\bR^3)=\mathcal{F}^{-1}(\mathcal{K}^{\alpha})$, and the existence and uniqueness of the solution to the Cauchy problem \eqref{1.1}--\eqref{1.3}--\eqref{1.4}  in the space $\tilde{P}_{\alpha}(\bR^3)$ follow from the argument in previous works \cite{morimoto2012remark,morimoto2015new}. For inverse power law potentials, the existence and uniqueness of the solution to the Cauchy problem \eqref{1.1}--\eqref{1.2}--\eqref{1.4} can be referred to \cite{morimoto2016measure,cannone2010infinite}.
    And it follows from definition 1.1 in \cite{morimoto2016measure}, here we can give the definition of the measure valued solution of the Cauchy problem \eqref{1.1}--\eqref{1.3}--\eqref{1.4}.
	\begin{defi}	\label{def1.1}
		Let $B(\cdot)$ satisfy \eqref{1.3} and $F_{0}\in P_{\alpha}(\bR^3)$ for $\alpha\in(0,2]$.  For $\psi\in C_b^2(\bR^3)$, set
		\[
		\Delta\psi=\psi(v')+\psi(v'_*)-\psi(v)-\psi(v_*).
		\]
		We say $F_{t}\in C([0,+\infty);P_{\alpha}(\bR^3))$ is a measure valued solution to the Cauchy problem \eqref{1.1}--\eqref{1.4} if it satisfies
		\begin{enumerate}
		    \item For every $\psi(v)\in C_{b}^{2}(\bR^3)$ and $t>0$,
		          $$
		          \begin{aligned}
			         \int_{0}^{t}\iint_{\bR^3\times \bR^3}\int_{\bS^2}B(|v-v_*|,\sigma)\Delta\psi d\sigma dF_{\tau}(v)dF_{\tau}(v_*)d\tau \mbox{ is finite}.
		          \end{aligned}
		          $$
            \item $\forall\psi(v)\in C_{b}^{2}(\bR^3)$, we have
			\begin{equation}
			\int_{\bR^3} \psi dF_{t}-	\int_{\bR^3} \psi dF_{0}=\frac{1}{2}\int_0^{t}\iiint_{\bR^3\times \bR^3\times \bS^2}B(|v-v_*|,\sigma)\Delta\psi d\sigma dF_{\tau}(v)dF_{\tau}(v_*)d\tau,
			\label{1.5}
			\end{equation}
			 where $\Delta\psi$ is defined as above.
            \item If $\alpha>1$, then the momentum conservation law holds:
		\begin{align*}
		    \forall t\geq0,\  \int v_{j}dF_{t}(v)=\int v_{j}dF_{0}(v), \ j=1,2,3.
		\end{align*}
        \item If $\alpha=2$, then $F_{t}\in C([0,\infty);P_{2}(\bR^3))$ and the energy conservation law holds:
		\begin{align*}
			\forall t\geq0,\  \int |v|^2dF_{t}(v)=\int |v|^2dF_{0}(v).
		\end{align*}	
		\end{enumerate}
	\end{defi}
 Since the collision cross section $B(v-v_\ast,\sigma)$ satisfies \eqref{1.3}, which implies that for any $\alpha>0$, we have
 $$
 \int_{\bS^2}b(\cos\theta)\sin^{\alpha}\theta d\sigma<+\infty.
 $$
 Therefore, we can use the same analytical techniques as Theorem 2.5 in Reference \cite{cannone2010infinite} to obtain the existence of the Cauchy problem \eqref{1.1}--\eqref{1.3}--\eqref{1.4}.\par 
 \begin{prop}\label{pro1.1}
 Let $B(\cdot)$ satisfy \eqref{1.3} and $F_{0}\in P_{\alpha}(\bR^3)$ for $\alpha\in(0,2]$. Then there exists a unique solution $F_t \in 
C([0,+\infty); P_{\alpha}(\bR^3))$ to the Cauchy problem \eqref{1.1}--\eqref{1.4} with the initial datum $F_0$. And for any $T>0$, there exists a $C_T>0$ such that
\begin{align}\label{alpha-moment-estimate}
\sup_{t\in[0,T]}\int|v|^{\alpha}dF_t(v)\leq C_T\int|v|^{\alpha}dF_0(v).
\end{align}
 \end{prop}
We remark that the continuity in time stated in Proposition \ref{pro1.1} is in the probability sense, that is, for any test function $\psi(v)\in C_b(\mathbb{R}^3)$, the mapping 
$$t\mapsto \int_{\bR^3} \psi(v) dF_t$$
is continuous.

 The main purpose of this paper is to prove the Villani conjecture for the Debye--Yukawa potential. The conjecture asserts that for any initial datum that is not a single Dirac mass, the corresponding measure valued solution becomes infinitely smooth in the Sobolev sense for all positive times. This conjecture has been proved in \cite{morimoto2015smoothing} for inverse power law potentials. In \cite{wang2024existence}, the smoothing effect of measure valued solutions was established under certain additional conditions on the initial datum, but the Villani conjecture was not proved there. In this paper, we prove this conjecture for the Debye--Yukawa potential. The main result is the following.
\begin{theo}\label{thm1.2}
	If the initial datum $F_{0}\in P_{\alpha}(\bR^3)$ for $\alpha\in(0,2]$ is not a single Dirac mass, and $F(t,v)\in C([0,\infty);P_{\alpha}(\bR^3))$ is the solution to the Cauchy problem \eqref{1.1}--\eqref{1.3}--\eqref{1.4}, then for any $t>0$, $F(t,\cdot)\in H^{\infty}(\bR^3)$.
\end{theo}
\par

\begin{rema}
(i) For any $t>0$, the support of $F(t, v)$ is $\bR^3$ if the initial datum $F_0$ is not a single Dirac mass. For $\alpha=2$, the proof is similar to the argument of Fournier \cite{fournier2015finiteness}; for $\alpha\in(0,2)$, we use the reverse Povzner inequality of Section 2.
\par
(ii) The key step in the proof of the smoothing effect is the coercivity estimate \cite{morimoto2015smoothing,wang2024existence}; a new proof of it is given in Section 3.
\end{rema}

The rest of the paper is organized as follows. In Section 2, we establish a reverse Povzner inequality which plays a crucial role in ruling out the possibility that the solution becomes a Dirac mass at any positive time, even if the initial datum has only finite fractional moments. Section 3 is devoted to the analysis of the coercivity estimate in Fourier space, which captures the regularization mechanism induced by angular non-cutoff interactions. Based on these tools, we prove the main result, Theorem \ref{thm1.2}, on the instantaneous smoothing effect of measure-valued solutions.

	\section{Reverse Povzner Inequality}
	 In this section, we establish a reverse Povzner inequality, which provides a lower bound for the energy at positive times. This estimate is essential for proving that the solution cannot concentrate into a Dirac mass if the initial datum is not a single Dirac mass, even when only fractional moments are finite.
	
	\begin{prop}\label{prop2.1}
		Assume $F_{0}\in P_{\alpha}(\bR^3)$, $1<\alpha<2$, $F(t,v)\in C([0,\infty);P_{\alpha}(\bR^3))$ is the solution to the Cauchy problem \eqref{1.1}--\eqref{1.3}--\eqref{1.4} corresponding to $F_0$ and there exists $t_0>0$ such that
		$$
		\int |v|^2dF_{t_0}(v)<+\infty.
		$$   
        Then
		\begin{equation}\label{2.1}
			\int|v|^2dF_{t_0}(v)\geq \int|v|^2dF_{0}(v).	
		\end{equation}	
	\end{prop}
	\begin{proof}
		For any $\delta>0$, let $\Psi(x)=\frac{x}{1+\delta x}$ for $x>0$, and then it follows from the equation that
			\begin{align*}
				&\frac{d}{dt}\int\Psi(|v|^2)dF_{t}(v)\\
                &=\frac{1}{2}\iiint_{\bR^3\times \bR^3\times \bS^2}b(\cos\theta)[\Psi(|v'|^2)+\Psi(|v'_\ast|^2)-\Psi(|v|^2)-\Psi(|v_*|^2)]d\sigma dF_{t}(v)dF_{t}(v_*)\\
                &=\frac{1}{2}\iiint_{\bR^3\times \bR^3\times \bS^2}b(\cos\theta)[G(v,v_*)-H(v,v_*)]d\sigma dF_{t}(v)dF_{t}(v_*),
			\end{align*}
		where
    \begin{align*}
			G(v,v_*)&=\Psi(|v'|^2)+\Psi(|v'_\ast|^2)-\Psi(|v|^2\sin^2\frac{\theta}{2}+|v_*|^2\cos^2\frac{\theta}{2})-\Psi(|v_*|^2\sin^2\frac{\theta}{2}+|v|^2\cos^2\frac{\theta}{2}),\\
            H(v,v_*)&=\Psi(|v|^2)+\Psi(|v_\ast|^2)-\Psi(|v|^2\sin^2\frac{\theta}{2}+|v_*|^2\cos^2\frac{\theta}{2})-\Psi(|v_*|^2\sin^2\frac{\theta}{2}+|v|^2\cos^2\frac{\theta}{2}).
	\end{align*}	

		Notice that
        \begin{equation}\label{2.2}
            \Psi''(x)=\frac{-2\delta}{(1+\delta x)^3}<0, \mbox{ and } |\Psi''(x)|\leq \frac{2\delta}{(1+\delta x)^{2-\alpha}},
        \end{equation}
		where $1<\alpha<2$. Then, because of the concavity of $\Psi(x)$, we have
		$$
		H(v,v_*)\leq 0.
		$$
		And by using the Povzner inequality studied in \cite{mischler1999spatially}, for $\sigma\in \bS^2$, it can be written as
		$$
\sigma=\mathbf{k}\cos\theta+\sin\theta(\mathbf{h}\cos\varphi+\mathbf{i}\sin\varphi),\ \theta\in[0,\frac{\pi}{2}],\ \varphi\in [-\pi,\pi),
		$$where
		$$
		\mathbf{k}=\frac{v-v_*}{|v-v_*|},\ 	\mathbf{i}=\frac{v\times v_*}{|v\times v_*|},\ \mathbf{h}=\mathbf{i}\times \mathbf{k}=\frac{((v-v_*)\cdot v)v_*-((v-v_*)\cdot v_*)v}{|v-v_*||v\times v_*|}.
		$$
			It follows from $(v+v_*)\perp\mathbf{i}$ and the definition of $\mathbf{h}$ that
		$$
		\begin{aligned}
			|v'|^2&=|\frac{v+v_*}{2}|^2+|\frac{v-v_*}{2}|^2+\frac{|v-v_*|}{2}(v+v_*)\cdot\sigma\\
			&=\frac{1}{4}(2|v|^2+2|v_*|^2)+\frac{|v-v_*|}{2}((v+v_*)\cdot(\cos\theta\mathbf{k}+\sin\theta \cos\varphi\mathbf{h}))\\&=\frac{1}{2}(|v|^2+|v_*|^2)+\frac{\cos\theta}{2}(|v|^2-|v_*|^2)\\&\ \ \ +\frac{\sin\theta \cos\varphi}{2|v\times v_*|}\{(v+v_*)\cdot(((v-v_*)\cdot v)v_*-((v-v_*)\cdot v_*)v)\}\\&=\frac{|v|^2(1+\cos\theta)}{2}+\frac{|v_*|^2(1-\cos\theta)}{2}+|v\times v_*|\sin\theta \cos\varphi.
		\end{aligned}
		$$
		Therefore, we have
		\begin{equation}
			|v'|^2=|v|^2\cos^2\frac{\theta}{2}+|v_*|^2\sin^2\frac{\theta}{2}+|v\times v_*|\sin\theta\cos\varphi \\
			:=Y(\theta)+Z(\theta)\cos\varphi,
			\label{2.3} 
		\end{equation}
        where 
        \begin{align*}
            Y(\theta)=|v|^2\cos^2\frac{\theta}{2}+|v_*|^2\sin^2\frac{\theta}{2},\ 
            Z(\theta)=|v\times v_*|\sin\theta.
        \end{align*}
			Similarly, 
            \begin{equation}
				|v'_{\ast}|^2=|v_*|^2\cos^2\frac{\theta}{2}+|v|^2\sin^2\frac{\theta}{2}-|v\times v_*|\sin\theta\cos\varphi:=Y(\pi-\theta)-Z(\theta)\cos\varphi.
                \label{2.4}
            \end{equation}
			Therefore,
			\begin{align*}
				&\int_{\bS^2}b(\cos\theta)G(v,v_*)d\sigma\\
                &=\int_{\bS^2}b(\cos\theta)[\Psi(Y(\theta)+Z(\theta)\cos\varphi)+\Psi(Y(\pi-\theta)-Z(\theta)\cos\varphi)-\Psi(Y(\theta))-\Psi(Y(\pi-\theta))]d\sigma\\
                &=2\int_{0}^{\frac{\pi}{2}}b(\cos\theta)\sin\theta\int_{0}^{\pi}\Psi(Y(\theta)+Z(\theta)\cos\varphi)+\Psi(Y(\pi-\theta)-Z(\theta)\cos\varphi)d\varphi d\theta\\
                &\ \ \ -2\pi\int_{0}^{\frac{\pi}{2}}b(\cos\theta)\sin\theta(\Psi(Y(\theta))+\Psi(Y(\pi-\theta)))d\theta.
			\end{align*}
The last equality follows from the spherical coordinate representation of the surface measure on $\mathbb{S}^2$ and the even parity of the integrand with respect to the azimuthal angle $\varphi$ over $[-\pi,\pi]$. Next, we notice that
	\begin{align*}
	    &\int_{0}^{\pi}\Psi(Y(\theta)+Z(\theta)\cos\varphi)d
		\varphi\\&=(\int_{0}^{\frac{\pi}{2}}+\int_{\frac{\pi}{2}}^{\pi})\Psi(Y(\theta)+Z(\theta)\cos\varphi)d
		\varphi\\&=\int_{0}^{\frac{\pi}{2}}\Psi(Y(\theta)+Z(\theta)\cos\varphi)+\Psi(Y(\theta)-Z(\theta)\cos\varphi)d\varphi.
	\end{align*}
Then, by integration by parts, we have
	$$
	\begin{aligned}
     &\int_{0}^{\pi}\Psi(Y(\theta)+Z(\theta)\cos\varphi)d
		\varphi\\
		&=[\varphi(\Psi(Y(\theta)+Z(\theta)\cos\varphi)+\Psi(Y(\theta)-Z(\theta)\cos\varphi))]_{0}^{\frac{\pi}{2}}\\&\ \ \ +Z(\theta)\int_{0}^{\frac{\pi}{2}}\varphi \sin\varphi[\Psi'(Y(\theta)+Z(\theta)\cos\varphi)-\Psi'(Y(\theta)-Z(\theta)\cos\varphi)]d\varphi\\&=\pi\Psi(Y(\theta))\\ &\ \ \ +Z^2(\theta)\int_{0}^{\frac{\pi}{2}}(\sin\varphi-\varphi\cos\varphi) \sin\varphi[\Psi''(Y(\theta)+Z(\theta)\cos\varphi)+\Psi''(Y(\theta)-Z(\theta)\cos\varphi)]d\varphi.
	\end{aligned}
	$$
	Similarly, we obtain 
	$$
	\begin{aligned}
		&\int_{0}^{\pi}\Psi(Y(\pi-\theta)-Z(\theta)\cos\varphi)d
		\varphi\\&=\pi\Psi(Y(\pi-\theta))\\ &\ \ \ +Z^2(\theta)\int_{0}^{\frac{\pi}{2}}(\sin\varphi-\varphi\cos\varphi) \sin\varphi[\Psi''(Y(\pi-\theta)+Z(\theta)\cos\varphi)+\Psi''(Y(\pi-\theta)-Z(\theta)\cos\varphi)]d\varphi.
	\end{aligned}
	$$
Therefore,
    $$
	\begin{aligned}
	&\int_{\bS^2}b(\cos\theta)G(v,v_*)d\sigma\\&=2\int_{0}^{\frac{\pi}{2}}b(\cos\theta)\sin\theta Z^2(\theta)\int_{0}^{\frac{\pi}{2}}(\sin\varphi-\varphi\cos\varphi) \sin\varphi[\Psi''(Y(\theta)+Z(\theta)\cos\varphi)+\Psi''(Y(\theta)-Z(\theta)\cos\varphi)\\&\ \ \ +
		\Psi''(Y(\pi-\theta)+Z(\theta)\cos\varphi)+\Psi''(Y(\pi-\theta)-Z(\theta)\cos\varphi)]d\varphi d\theta,
	\end{aligned}
	$$
    then, it follows from \eqref{2.2} that 
    \begin{align*}
    \int_{\bS^2}b(\cos\theta)G(v,v_*&)d\sigma\lesssim \delta \int_{0}^{\frac{\pi}{2}}b(\cos\theta)\sin\theta Z^2(\theta)\int_{0}^{\frac{\pi}{2}}(\sin\varphi-\varphi\cos\varphi) \sin\varphi \times\\&\ \ \ \ \ \ \ \ \ \ \ [\frac{1}{(1+\delta Y(\theta)+\delta Z(\theta)\cos\varphi)^{2-\alpha}}+\frac{1}{(1+\delta Y(\theta)-\delta Z(\theta)\cos\varphi)^{2-\alpha}}]d\theta d\varphi\\&\lesssim Z_{0}^{2-\alpha} \delta^{\alpha-1} \int_{0}^{\frac{\pi}{2}}b(\cos\theta)\sin\theta Z^{\alpha}(\theta)\int_{0}^{\frac{\pi}{2}}\varphi^4  [\frac{1}{(1+Z_{0}\cos\varphi)^{2-\alpha}}+\frac{1}{(1-Z_{0}\cos\varphi)^{2-\alpha}}]d\theta d\varphi\\&\lesssim \delta^{\alpha-1}|v|^{\alpha
		}|v_*|^{\alpha},
    \end{align*}
		where $Z_{0}=\frac{\delta Z(\theta)}{1+\delta Y(\theta)}\in [0,1]$, and the last inequality follows from $Z^{\alpha}(\theta)\lesssim|v|^{\alpha
		}|v_*|^{\alpha}\theta^{\alpha}$ and
		$$
		\int_{0}^{\frac{\pi}{2}}b(\cos\theta) \theta^{\alpha+1}d\theta<\infty.
		$$
		Hence, 
        \begin{align*}
				\frac{d}{dt}\int\Psi(|v|^2)dF_{t}(v)\geq-\frac{1}{2}|\iiint_{\bR^3\times \bR^3\times \bS^2}b(\cos\theta)G(v,v_*)d\sigma dF_{t}(v)dF_{t}(v_*)|,
			\end{align*}
        then, by  \eqref{alpha-moment-estimate} with $T=t_0$,
		$$	
		\begin{aligned}
			\dfrac{d}{dt}\int\Psi(|v|^2)dF_{t}(v)\geq -\delta^{\alpha-1}C(\int |v|^{\alpha}dF_{t}(v))^2\geq -\delta^{\alpha-1}CC_{t_0}^2,
		\end{aligned}
		$$
		holds for any $t\in(0,t_0)$, where $C$ is a constant independent of $\delta$. 
        Consequently,
        $$\int\Psi(|v|^2)dF_{t_0}(v)\ge \int\Psi(|v|^2)dF_{0}(v)-\delta^{\a-1}CC_{t_0}^2t_0.$$
        Since $1<\alpha<2$, taking the limit as $\delta\to0$, we can obtain the estimate \eqref{2.1}.
		
	\end{proof}

	\begin{prop}\label{prop2.2}
			Assume $F_{0}\in P_{\alpha}(\bR^3)$, $0<\alpha\leq 1$, $F(t,v)\in C([0,\infty);P_{\alpha}(\bR^3))$ is the solution to the Cauchy problem \eqref{1.1}--\eqref{1.3}--\eqref{1.4} corresponding to $F_0$ and there exists $t_0>0$ such that
		$$
		\int|v|^{2\alpha}dF_{t_0}(v)<+\infty.
		$$
        Then there exists $C_{t_0}'>0$ such that
		\begin{equation}
		\sup_{[0,t_0]}\int|v|^{2\alpha}dF_{t}(v)\leq C_{t_0}'.
			\label{2.5}
		\end{equation}
	\end{prop}
	\begin{proof}
	For any $\delta>0$, $x>0$, define $\Psi_{\alpha}(x)=\frac{x^{\alpha}}{1+\delta x^{\alpha}}$. Then for $0<\alpha\leq1$, we have
		\begin{equation}
			\begin{split}
				\Psi''_{\alpha}(x)=\frac{\alpha x^{\alpha-2}}{(1+\delta x^{\alpha})^3}&[(\alpha-1)-\delta (\alpha+1)x^{\alpha}]\leq 0,\\&
				|\Psi''_{\alpha}(x)|\leq C_{0}  x^{\alpha-2},
			\end{split}
			\label{2.6}
		\end{equation}
	where $C_{0}$ is a constant independent of $\delta$. Then following the proof of Proposition \ref{prop2.1}, we have
		$$
	\begin{aligned}
		&\frac{d}{dt}\int\Psi_{\alpha}(\langle v\rangle^2)dF_{t}(v)\\&=\frac{1}{2}\iiint_{\bR^3\times \bR^3\times \bS^2}b(\cos\theta)[\Psi_{\alpha}(\langle v'\rangle^2)+\Psi_{\alpha}(\langle v'_\ast\rangle^2)-\Psi_{\alpha}(\langle v\rangle^2)-\Psi_{\alpha}(\langle v_*\rangle^2)]d\sigma dF_{t}(v)dF_{t}(v_*)\\&=\frac{1}{2}\iiint_{\bR^3\times \bR^3\times \bS^2}b(\cos\theta)[G(v,v_*)-H(v,v_*)]d\sigma dF_{t}(v)dF_{t}(v_*),
	\end{aligned}
	$$
	where
	$$	
	\begin{aligned}
		&G(v,v_*)=\Psi_{\alpha}(\langle v'\rangle^2)+\Psi_{\alpha}(\langle v'_\ast\rangle^2)-\Psi_{\alpha}(\langle v\rangle^2\sin^2\frac{\theta}{2}+\langle v_*\rangle^2\cos^2\frac{\theta}{2})-\Psi_{\alpha}(\langle v_*\rangle^2\sin^2\frac{\theta}{2}+\langle v\rangle^2\cos^2\frac{\theta}{2}),\\&H(v,v_*)=\Psi_{\alpha}(\langle v\rangle^2)+\Psi_{\alpha}(\langle v_\ast\rangle^2)-\Psi_{\alpha}(\langle v\rangle^2\sin^2\frac{\theta}{2}+\langle v_*\rangle^2\cos^2\frac{\theta}{2})-\Psi_{\alpha}(\langle v_*\rangle^2\sin^2\frac{\theta}{2}+\langle v\rangle^2\cos^2\frac{\theta}{2}).
	\end{aligned}	
	$$
		It follows from \eqref{2.3} and \eqref{2.4} that
        \begin{align*}
            \langle v'\rangle^2&=1+|v'|^2\\&=1+|v|^2\cos^2\frac{\theta}{2}+|v_*|^2\sin^2\frac{\theta}{2}+|v\times v_*|\sin\theta\cos\varphi\\&=\langle v\rangle^2\cos^2\frac{\theta}{2}+\langle v_*\rangle^2\sin^2\frac{\theta}{2}+|v\times v_*|\sin\theta\cos\varphi\\&:=Y_1(\theta)+Z(\theta)\cos\varphi,\\
            \langle v'_{\ast}\rangle^2&=\langle v_*\rangle^2\cos^2\frac{\theta}{2}+\langle v\rangle^2\sin^2\frac{\theta}{2}-|v\times v_*|\sin\theta\cos\varphi\\&:=Y_1(\pi-\theta)-Z(\theta)\cos\varphi,
        \end{align*}
	where $Y_1(\theta)=1+Y(\theta)$, $Y_1(\pi-\theta)=1+Y(\pi-\theta)$. Due to the concavity of $\Psi_{\alpha}(x)$, we can also obtain$$
	H(v,v_*)\leq 0.
	$$\par 
	Then, following the proof of Propositions \ref{prop2.1} and \eqref{2.6}, we have
	\begin{align*}
		&\int_{\bS^2}b(\cos\theta)G(v,v_*)d\sigma\\
        &=2\int_{0}^{\frac{\pi}{2}}b(\cos\theta)\sin\theta Z^2(\theta)\int_{0}^{\frac{\pi}{2}}(\sin\varphi-\varphi\cos\varphi) \sin\varphi[\Psi_{\alpha}''(Y_1(\theta)+Z(\theta)\cos\varphi)+\Psi_{\alpha}''(Y_1(\theta)-Z(\theta)\cos\varphi)\\
        &\ \ \ +
		\Psi_{\alpha}''(Y_1(\pi-\theta)+Z(\theta)\cos\varphi)+\Psi_{\alpha}''(Y_1(\pi-\theta)-Z(\theta)\cos\varphi)]d\varphi d\theta.\\
        &\lesssim \int_{0}^{\frac{\pi}{2}}b(\cos\theta)\sin\theta Z^{\alpha}(\theta)Z_{0}^{2-\alpha}\int_{0}^{\frac{\pi}{2}}\varphi^4 [(1+Z_{0}\cos\varphi)^{\alpha-2}+(1-Z_{0}\cos\varphi)^{\alpha-2}]d\varphi d\theta\\&\leq C_{1}|v|^{\alpha}|v_*|^{\alpha},
	\end{align*}
where $C_{1}$ is a constant independent of $\delta$, and $$Z_{0}=\frac{Z(\theta)}{Y_1(\theta)}=\frac{|v\times v_*|\sin\theta}{1+|v|^2\cos^2\frac{\theta}{2}+|v_*|^2 \sin^2\frac{\theta}{2}}\in[0,1].$$
In fact, it follows from $0<\alpha\leq1$ that
$$\int_{0}^{\frac{\pi}{2}}\varphi^4 [(1+Z_{0}\cos\varphi)^{\alpha-2}+(1-Z_{0}\cos\varphi)^{\alpha-2}]d\varphi$$ is integrable.
	 Therefore, for any $t\in(0,t_0)$,
	$$
	\begin{aligned}
		&\int\Psi_{\alpha}(\langle v\rangle^2)dF_{t_0}(v)\\
        &\geq\int\Psi_{\alpha}(\langle v\rangle^2)dF_{t}(v)-\frac{1}{2}|\int_t^{t_0}d\tau\iiint_{\bR^3\times \bR^3\times \bS^2}b(\cos\theta)G(v,v_*)d\sigma dF_{\tau}(v)dF_{\tau}(v_*)|\\
        &\geq\int\Psi_{\alpha}(\langle v\rangle^2)dF_{t}(v)-C_{1}\int_t^{t_0}d\tau(\int_{ \bR^3}|v|^{\alpha} dF_{\tau}(v))^2\\
        &\geq\int\Psi_{\alpha}(\langle v\rangle^2)dF_{t}(v)-C_1C^2_{t_0}t_0,
	\end{aligned}
	$$
	where $C_1,C_{t_0}$ are independent of $\delta$. Letting $\delta\to0$, we obtain \eqref{2.5}, if one takes
    $$C_{t_0}'=\int\langle v\rangle^{2\a}dF_{t_0}+C_1C_{t_0}^2 t_0.$$
	
	\end{proof}

	\section{Smoothing Effect of the Measure Valued Solution}
	 In this section, we give a new insight into the coercivity estimate, which is the key step to prove the smoothing effect of the measure valued solution to the Cauchy problem \eqref{1.1}--\eqref{1.3}--\eqref{1.4}. We first prove that the support of the solution is the whole space at any positive time, and then establish the coercivity estimate in Fourier space.
	\begin{prop}\label{prop3.1}
		If $F_{0}\in P_{\alpha}(\bR^3)$ for $\alpha\in(0,2]$ is not a single Dirac mass, then for any weak solution $F_{t}\in C([0,\infty);P_{\alpha}(\bR^3))$ in the sense of Definition \ref{def1.1} to the Cauchy problem \eqref{1.1}--\eqref{1.3}--\eqref{1.4}, we have 
		\begin{align*}
			\mathrm{supp} \ F_{t}=\bR^3,\mbox{ for all }t>0.
		\end{align*}	
	\end{prop}
	\begin{proof}
		Following the proof of Theorem 1.2 in \cite{fournier2015finiteness}, we will divide the proof into several steps.
		\par 
	Firstly, we should claim that for all $\psi(v)\in C(\bR^3)$ satisfies the growth condition $ \dis\sup_{v\in \bR^3}\frac{|\psi(v)|}{\langle v\rangle^{\alpha}}<\infty$, then
	\begin{align*}
	t\to\int_{\bR^3}\psi(v)dF_{t}(v)\in C^1([0,+\infty)).
	\end{align*}
	It follows from \eqref{1.5} that we just need to prove the map
	$$
	\begin{aligned}
		t\to \iiint_{\bR^3\times \bR^3\times \bS^2}b(\cos\theta)\Delta\psi d\sigma dF_{t}(v)dF_{t}(v_*)
	\end{aligned}
	$$
	belongs to $C([0,+\infty))$.  In fact, since $F_{t}\in C([0,\infty);P_{\alpha}(\bR^3))$, for any $\varepsilon>0$ there exists $R > 0$ independent
	of $t$ such that
	$$
	\begin{aligned}
		|\iint_{\{|v|\geq R\}\bigcup\{|v_*|\geq R\}}\int_{\bS^2}b(\cos\theta)\Delta\psi d\sigma dF_{t}(v)dF_{t}(v_*)|<\varepsilon,
	\end{aligned}
	$$
	while on the compact 
	$$
	K_R = \{(v, v_*) :|v|\leq R, |v_*|\leq R\},
	$$
	we can use Weierstrass approximation argument (see Lemma 3.2 in \cite{morimoto2016measure}) to obtain the continuity.
	\par 
	Next, we claim that if $F_{0}\in P_{\alpha}(\bR^3)$, $ 0<\alpha\leq 2$ is not a single Dirac mass, then for all $t > 0$, $F_{t}$ is also not a Dirac mass. We discuss it in three cases.
    
    Case 1: For $\alpha=2$, it is immediate from the conservation of momentum and energy in Definition \ref{def1.1} and the fact that $F_{0}$ is not a Dirac mass: 
	$$
	\begin{aligned}
		\int_{\bR^3}|v-v_{0}|^2 dF_{t}(v)=\int_{\bR^3}|v-v_{0}|^2 dF_{0}(v)>0,\ \forall v_{0}\in \bR^3,\ \forall t>0,
	\end{aligned}
	$$
	that is, $F_{t}$ is also not a Dirac mass. 
    
    Case 2: For $1<\a<2$,  assume that there exists $t_0>0$, $v_{0}\in \bR^3$, $F_{t_0}=\delta_{v_0}$. Then, by using Proposition \ref{prop2.1}, we have
	$$
	\begin{aligned}
		 \int|v|^2dF_{0}(v)\leq\int|v|^2dF_{t_0}(v)=|v_{0}|^2<\infty.
	\end{aligned}
	$$
    Then $F_0\in P_2$, and we return to the Case 1.
    
    Case 3: For $0<\alpha\leq 1$, assume that there exists $t_0>0$, $v_{0}\in \bR^3$, $F_{t_0}=\delta_{v_0}$. Then it follows from Proposition \ref{prop2.2}, we have \eqref{2.5}
	$$
	\sup_{[0,t_0]}\int|v|^{2\alpha}dF_{t}(v)\leq C'_{t_0}.
	$$
	If $2\alpha\le1$, repeating the proof process of Proposition \ref{prop2.2}, we obtain
	$$
	\sup_{[0,t_0]}\int|v|^{2^2\alpha}dF_{t}(v)\leq C''_{t_0}.
	$$
	After repeating $n_0$ times, with $n_0$ satisfying $1<2^{n_{0}}\alpha\leq 2$, it returns to the argument in Case 2.
    \par 
	The remaining steps are the same as the proof of \cite[Theorem 1.2]{fournier2015finiteness}.

	\end{proof}

Now, we can prove the following coercivity estimate, which is a key to the proof of Theorem \ref{thm1.2}. By the Bobylev formula\cite{alexandre2000entropy}, the Cauchy problem \eqref{1.1}--\eqref{1.3}--\eqref{1.4} is reduced to 
\begin{equation}
	\begin{cases}
		\partial_{t} \varphi(t,\xi)=\int_{\bS^2} b(\frac{\xi\cdot \sigma}{|\xi|})(\varphi(t,\xi^+)\varphi(t,\xi^-)-\varphi(t,\xi)\varphi(t,0)) d\sigma\\
		\varphi(0,\xi)=	\int_{\bR^3} e^{-iv \cdot \xi}dF_{0}(v)
	\end{cases},
	\label{3.1}
\end{equation}
where $\varphi(t,\xi)= \mathcal{F}(F_{t})(\xi)$ and 
\[
\xi^+=\frac{\xi+|\xi|\sigma}{2},\qquad
\xi^-=\frac{\xi-|\xi|\sigma}{2}.
\]

Denote the regions
\[
\Omega_i=\left\{x\in\bR^3:\ \angle(x,e_i)\le \frac{\pi}{8},\ |x|\ge1\right\},
\qquad i=1,2,
\]
and 
\[
\Omega_i'=\left\{x\in\bR^3:\ \angle(x,e_i)\le \frac{\pi}{16},\ |x|\ge1\right\},
\qquad i=1,2.
\]
Since $\Omega_1 \cap \Omega_2 = \emptyset$, at least one of $\xi \notin \Omega_1$ and $\xi \notin \Omega_2$ must hold. The following lemma provides a key angular integral estimate, which will be used in the proof of the coercivity estimate.
\begin{lemm}[Angular integral]\label{lem:angular}
Let $|\xi|>1$. For $i=1,2$, if $\xi\notin\Omega_i$, then
\[
\int_{\bS^2} b\!\left(\frac{\xi\cdot\sigma}{|\xi|}\right)
\sin^2\frac{\xi^-\cdot(v-u)}2
\,d\sigma
\gtrsim
\mathbbm{1}_{\{(v-u)\in \Omega_i'\}}
\int_{1/|\xi|}^{\pi/2} b(\cos\theta)\sin\theta\,d\theta .
\]
\end{lemm}

\begin{proof}
Notice for any $\varepsilon>0$ small enough, 
\[
\int_{\bS^2} b\!\left(\frac{\xi\cdot\sigma}{|\xi|}\right)
\sin^2\frac{\xi^-\cdot(v-u)}2
\,d\sigma
\gtrsim
\int_{\bS^2} b\!\left(\frac{\xi\cdot\sigma}{|\xi|}\right)
\mathbbm{1}_{\{(v-u)\cdot\xi^-\notin
\bigcup_{n\in\bZ}[2n\pi-2\varepsilon,2n\pi+2\varepsilon]\}}
\,d\sigma.
\]
Choose an orthonormal basis
$\{\boldsymbol h_1,\boldsymbol h_2,\boldsymbol h_3=\xi/|\xi|\}$ and write
\[
\sigma=\sin\theta\cos\phi\,\boldsymbol h_1
+\sin\theta\sin\phi\,\boldsymbol h_2
+\cos\theta\,\boldsymbol h_3,
\]
\[
v-u=|v-u|\bigl(\sin\theta_1\cos\phi_1\,\boldsymbol h_1
+\sin\theta_1\sin\phi_1\,\boldsymbol h_2
+\cos\theta_1\,\boldsymbol h_3\bigr).
\]
Then
\[
(v-u)\cdot\sigma
=|v-u|\sin\theta\sin\theta_1\cos(\phi-\phi_1)
+\frac{(v-u)\cdot\xi}{|\xi|}\cos\theta.
\]
Since $\xi^-=(\xi-|\xi|\sigma)/2$, we have
\[
(v-u)\cdot\xi^-
=\frac12\Bigl((v-u)\cdot\xi-|\xi|(v-u)\cdot\sigma\Bigr).
\]
Substituting the above expression for $(v-u)\cdot\sigma$ gives
\[
(v-u)\cdot\xi^-
=\frac12\Bigl((v-u)\cdot\xi(1-\cos\theta)
-|\xi||v-u|\sin\theta\sin\theta_1\cos(\phi-\phi_1)\Bigr).
\]
Set
\[
L:=|\xi|\sin\theta|v-u|\sin\theta_1.
\]
We claim that if $v-u\in\Omega_i'$ and $\xi\notin\Omega_i$, then
\[
\sin\theta_1\ge\sin(\pi/16)=:C_0>0.
\]
Indeed, $v-u$ makes an angle at most $\pi/16$ with $e_i$, while $\xi$ makes an angle larger than $\pi/8$ with $e_i$. Hence the angle between $v-u$ and $\xi$ is at least $\pi/8-\pi/16=\pi/16$, which is exactly the angle $\theta_1$. Therefore, since $|\xi|\sin\theta\gtrsim1$ and $|v-u|\ge1$, we get
\[
L\ge C_0'>0.
\]
For fixed $\theta$, define the admissible set
\[
\cG_\theta
:=\left\{\phi\in[0,2\pi]:
(v-u)\cdot\xi^-\notin
\bigcup_{n\in\bZ}[2n\pi-2\varepsilon,2n\pi+2\varepsilon]
\right\}.
\]
The condition $(v-u)\cdot\xi^-\notin
\bigcup_{n\in\bZ}[2n\pi-2\varepsilon,2n\pi+2\varepsilon]$
is equivalent to
\[
\cos(\phi-\phi_1)\notin
\bigcup_{n\in\bZ}[A_n-\delta_n,A_n+\delta_n],
\]
where
\[
A_n=
\frac{4n\pi+(v-u)\cdot\xi(1-\cos\theta)}{L},
\qquad
\delta_n=
\frac{4\varepsilon}{L}.
\]
The centers $A_n$ are spaced by $4\pi/L$, while each forbidden
interval has length $8\varepsilon/L$. Hence the total length of all
forbidden intervals in the $\cos$-variable is at most
\[
\sum_n \frac{8\varepsilon}{L}
\lesssim \varepsilon,
\]
because the number of $n$ for which $A_n\in[-1,1]$ is $O(L)$.
Let \(E\subset[-1,1]\) be the union of these forbidden intervals, so that \(|E|\lesssim \varepsilon\).
After the change of variable \(x=\cos(\phi-\phi_1)\), the preimage of \(E\) in the \(\phi\)-variable has measure at most
\[
C_1 |E| + C_2 \sqrt{|E|},
\]
where \(C_1,C_2>0\) are absolute constants. 
Indeed, splitting \(E\) into \(E\cap\{|x|\le 1/2\}\) and \(E\cap\{|x|\ge 1/2\}\), the first part contributes at most \(O(|E|)\), while the second part, being contained in two intervals near \(\pm1\), contributes at most \(O(\sqrt{|E|})\).
Therefore, for \(\varepsilon>0\) sufficiently small, this measure is less than \(\pi/2\). Hence the admissible set \(\cG_\theta\) satisfies
\[
|\cG_\theta|\ge \frac{\pi}{2}.
\]
Consequently,
\[
\begin{aligned}
\int_{\bS^2} b\!\left(\frac{\xi\cdot\sigma}{|\xi|}\right)
\mathbbm{1}_{\cG_\theta}(\phi)\,d\sigma
&\ge
\int_{1/|\xi|}^{\pi/2} \int_{\cG_\theta}
b(\cos\theta)\sin\theta\,d\phi\,d\theta\\
&\gtrsim
\int_{1/|\xi|}^{\pi/2} b(\cos\theta)\sin\theta\,d\theta,
\end{aligned}
\]
where we used $|\cG_\theta|\ge\pi/2$ in the last step.
\end{proof}
\begin{prop}[coercivity estimate]\label{prop3.2}
Assume that $b(\cdot)$ is given in \eqref{1.3} with $m=2/s-1>0$, and $F_{0}\in P_{\alpha}(\bR^3)$ for $0< \alpha\leq 2$ is not a single Dirac mass, $\varphi(t,\xi)$ is the solution of the Cauchy problem \eqref{3.1}, then there exists $q(t)>0$ such that
\begin{align*}
	\int_{\bR^3}\int_{\bS^2} b(\frac{\xi\cdot \sigma}{|\xi|})(1-|\varphi(t,\xi^-)|)d\sigma|h(\xi)|^2d\xi+\int_{\bR^3}|h(\xi)|^2d\xi\geq q(t)\int_{\bR^3}(\log\langle\xi\rangle)^{m+1}|h(\xi)|^2d\xi
\end{align*}
holds for any $t>0$ and for any $h(\xi)\in L_{\alpha}^{2}(\bR^3)$, $h\neq0$.
\end{prop}

\begin{proof}
Set
\[
\xi^+=\frac{\xi+|\xi|\sigma}{2},\qquad
\xi^-=\frac{\xi-|\xi|\sigma}{2}.
\]
The proof is divided into four steps.

First, for any \(\xi\neq0\), since \(|\varphi(t,\xi)|\le1\),
\[
\begin{aligned}
1-|\varphi(t,\xi)|
&\ge \frac{1-|\varphi(t,\xi)|^2}{2}\\
&=\frac12\iint_{\bR^3\times\bR^3}
\bigl(1-\cos(\xi\cdot(v-u))\bigr)\,dF_t(u)dF_t(v)\\
&=\iint_{\bR^3\times\bR^3}
\sin^2\frac{\xi\cdot(v-u)}2\,dF_t(u)dF_t(v).
\end{aligned}
\]
Hence
\begin{equation}\label{3.0}
1-|\varphi(t,\xi^-)|
\ge \iint_{\bR^3\times\bR^3}
\sin^2\frac{\xi^-\cdot(v-u)}2\,dF_t(u)dF_t(v).
\end{equation}

Second, by Lemma \ref{lem:angular}, for $|\xi|>1$,
\begin{equation}\label{3.2}
\int_{\bS^2} b\!\left(\frac{\xi\cdot\sigma}{|\xi|}\right)
\bigl(1-|\varphi(t,\xi^-)|\bigr)\,d\sigma
\gtrsim
c_1(t)(\log\langle\xi\rangle)^{m+1},
\end{equation}
where
\[
c_1(t):=\min_{i=1,2}\iint_{v-u\in\Omega_i'}dF_t(u)dF_t(v).
\]
Indeed, combining Lemma \ref{lem:angular} with \eqref{1.3} and
$m=2/s-1>0$, we have
\[
\int_{1/|\xi|}^{\pi/2} b(\cos\theta)\sin\theta\,d\theta
\gtrsim (\log\langle\xi\rangle)^{m+1}.
\]

Third, by Proposition \ref{prop3.1}, $\operatorname{supp}F_t=\bR^3$ for every $t>0$. We claim that the constant $c_1(t)$ introduced in \eqref{3.2} is strictly positive, that is,
\[
c_1(t)=\min_{i=1,2}\iint_{v-u\in\Omega_i'}dF_t(u)dF_t(v)>0.
\]
Indeed, fix any $x_0\in\bR^3$ and $r,R>0$. Choose $d>0$ such that
\[
d\ge R+r+1,\qquad d\sin\frac{\pi}{16}\ge R+r.
\]
The first condition ensures that the balls $B_r(x_0)$ and $B_R(x_0+de_i)$ are disjoint, while the second one guarantees that for all $u\in B_r(x_0)$ and $v\in B_R(x_0+de_i)$,
\[
v-u\in\Omega_i'.
\]
Hence
\[
B_r(x_0)\times B_R(x_0+de_i)
\subset \{(u,v)\in\bR^6:\ v-u\in\Omega_i'\}.
\]
Therefore
\[
\iint_{v-u\in\Omega_i'}dF_t(u)dF_t(v)
\ge F_t(B_r(x_0))\,F_t(B_R(x_0+de_i))>0,
\]
where the strict positivity follows from $\operatorname{supp}F_t=\bR^3$. Set
\[
c(t):=\min_{i=1,2}F_t(B_r(x_0))\,F_t(B_R(x_0+de_i))>0.
\]
The computation above shows $c_1(t)\ge c(t)$. Inserting this into \eqref{3.2}, we obtain, for $|\xi|>1$,
\begin{equation}\label{3.3}
\int_{\bS^2} b\!\left(\frac{\xi\cdot\sigma}{|\xi|}\right)
\bigl(1-|\varphi(t,\xi^-)|\bigr)\,d\sigma
\gtrsim
c_1(t)(\log\langle\xi\rangle)^{m+1}
\ge
c(t)(\log\langle\xi\rangle)^{m+1}.
\end{equation}

Finally, multiplying \eqref{3.3} by $|h(\xi)|^2$ and integrating over $|\xi|>1$, we obtain
\[
\int_{|\xi|>1}\int_{\bS^2} b\!\left(\frac{\xi\cdot\sigma}{|\xi|}\right)
\bigl(1-|\varphi(t,\xi^-)|\bigr)d\sigma\,|h(\xi)|^2d\xi
\gtrsim
c(t)\int_{|\xi|>1}(\log\langle\xi\rangle)^{m+1}|h(\xi)|^2d\xi.
\]
For $|\xi|\le 1$, since $\log\langle\xi\rangle\le \log 2$, there exists a constant $C_1>0$ such that
\[
(\log\langle\xi\rangle)^{m+1}\le C_1,\qquad |\xi|\le 1.
\]
Hence
\[
\int_{|\xi|\le 1}(\log\langle\xi\rangle)^{m+1}|h(\xi)|^2d\xi
\le C_1\int_{|\xi|\le 1}|h(\xi)|^2d\xi
\le C_1\int_{\bR^3}|h(\xi)|^2d\xi.
\]
Combining the two estimates, we get
\[
\begin{aligned}
&\int_{\bR^3}\int_{\bS^2} b\!\left(\frac{\xi\cdot\sigma}{|\xi|}\right)
\bigl(1-|\varphi(t,\xi^-)|\bigr)d\sigma\,|h(\xi)|^2d\xi
+\int_{\bR^3}|h(\xi)|^2d\xi\\
&\qquad\gtrsim
\min\{c(t),C_1^{-1}\}
\int_{\bR^3}(\log\langle\xi\rangle)^{m+1}|h(\xi)|^2d\xi.
\end{aligned}
\]
Set
\[
q(t):=\min\{c(t),C_1^{-1}\}>0.
\]
Then $q(t)>0$ and the desired coercivity estimate follows. This completes the proof.

\end{proof}
Now, we are ready to prove Theorem \ref{thm1.2}. Before proving Theorem \ref{thm1.2}, we clarify the regularity and integrability of $q(t)$, which are needed in the weighted energy estimate. 
Following the argument in \cite{wang2024existence}, the function $q(t)$ is positive for every $t>0$, measurable, and locally integrable. If necessary, we replace $q(t)$ by $\min\{q(t), e^{m+1}+1\}$ to ensure local boundedness. Therefore, we have
\[
\frac{d}{dt}\int_0^t q(\tau)d\tau = q(t), \quad \text{a.e. } t \in (0,T_0).
\]
This justifies the use of the Gronwall inequality in the subsequent proof.
\begin{proof}[Proof of Theorem \ref{thm1.2}] For $t \in [0, T_0].$ Define a time dependent weight function
	$$
M_{\delta}	(t,\xi)=\langle\xi\rangle^{N\int_0^tq(\tau)d\tau-4}(1+\delta|\xi|^2)^{-N_0}
,$$
where  $N_{0}=NT_{0}\sup_{t\in[0,T_{0}]}|q(t)|+2$, $N\in \mathbb{N}$, $\delta>0$. We multiply  $M_{\delta}^2	(t,\xi)\overline{\varphi(t,\xi)}$ by the equation \eqref{3.1} and integrate with respect to $\xi$ 
over $\bR^3$. Denote $\varphi^{\pm}=\varphi(t,\xi^{\pm})$, $M^+=M_{\delta}	(t,\xi^+)$ to simplify the notation and note
that
$$
\begin{aligned}
	&\partial_t \int_{\bR^3} |M_{\delta}\varphi|^2d\xi-2Nq(t)\int_{\bR^3} \log\langle\xi\rangle |M_{\delta}\varphi|^2d\xi =2\int_{\bR^3}\int_{\bS^2}b(\frac{\xi\cdot \sigma}{|\xi|})Re[(\varphi^{+}\varphi^{-} -\varphi)M_\delta^2\overline{\varphi}] d\sigma d\xi,
\end{aligned}
$$
	and 
		$$
	\begin{aligned}
		&-2Re[(\varphi^{+}\varphi^{-} -\varphi)M_\delta^2\overline{\varphi}]\\&=(|M_\delta \varphi|^2-2Re\{\varphi^{-}\varphi^{+}M^{+}\overline{M_\delta\varphi}\}+|M^{+} \varphi^{+}|^2)\\&\ \ \ \ +(|M_\delta \varphi|^2-|M^{+} \varphi^{+}|^2)\\
		&\ \ \ \ -2Re\{\varphi^{-}((M_\delta-M^{+})\varphi^{+})\overline{M_\delta\varphi}\}\\&=J_1+J_2+J_3.
	\end{aligned}
	$$
	Using the Cauchy--Schwarz inequality for the third term of $J_{1}$, we have
	$$
	\begin{aligned}
		J_1\geqslant(1-|\varphi^{-}|)(|M_\delta \varphi|^2+|M^{+} \varphi^{+}|^2)\geqslant(1-|\varphi^{-}|)|M_\delta \varphi|^2.
	\end{aligned}
	$$
	Therefore, by means of Proposition \ref{prop3.2}, we get
	$$
	\begin{aligned}
		\int_{\bR^3}\int_{\bS^2}b(\frac{\xi\cdot \sigma}{|\xi|})J_1 d\sigma d\xi+\int_{\bR^3} |M_\delta \varphi|^2d\xi
		\gtrsim q(t) \int_{\bR^3}(\log\langle\xi\rangle)^{m+1}|M_\delta \varphi|^2d\xi.
	\end{aligned}
	$$
	If we use the change of variable $\xi \rightarrow \xi 
	^+$ for the term $M^+\varphi^ 
	+$ in $J_{2}$, by the cancellation
	lemma (Lemma 1 in \cite{alexandre2000entropy}), we have
	$$
	\begin{aligned}
		|\int_{\bR^3}\int_{\bS^2}b(\frac{\xi\cdot \sigma}{|\xi|})J_2 d\sigma d\xi|&=2\pi|\int_{\bR^3}|M_\delta \varphi|^2(\int_0^\frac{\pi}{2}b(\cos\theta)\sin\theta(1-\frac{1}{\cos^3\frac{\theta}{2}})d\theta)d\xi|\\ &\lesssim\int_{\bR^3}|M_\delta \varphi|^2d\xi.
	\end{aligned}
	$$
	Since $|M_\delta-M^+|\lesssim \sin^2(\theta/2) M^+$ (see (3.4) in \cite{morimoto2009regularity}),  by the Cauchy--Schwarz inequality
	and the change of
	variable $\xi \rightarrow \xi^+$ for the term including $M^+\varphi^+$, we can also obtain
	$$
	\begin{aligned}
		|\int_{\bR^3}\int_{\bS^2}b(\frac{\xi\cdot \sigma}{|\xi|})J_3 d\sigma d\xi|&\lesssim \int_{\bR^3}\int_{\bS^2}b(\frac{\xi\cdot \sigma}{|\xi|})|M^+\varphi^+||M\varphi|\sin^2\frac{\theta}{2}d\sigma d\xi \\ 
		&\lesssim \int_{\bR^3}|M_\delta \varphi|^2d\xi.
	\end{aligned}
	$$
	 Since $m > 0$, $(\log\langle \xi \rangle )^m \rightarrow \infty$ as $| \xi |\rightarrow \infty$, we have 
	 $$
	 \begin{aligned}
	 	\partial_t \int_{\bR^3} |M_{\delta}\varphi|^2d\xi\lesssim\int_{\bR^3} |M_{\delta}\varphi|^2d\xi,
	 \end{aligned}
	 $$
and it follows from the Gronwall inequality, we have
$$
\begin{aligned}
	\int_{\bR^3}|\langle\xi\rangle^{N\int_0^tq(\tau)d\tau-4}(1+\delta|\xi|^2)^{-N_0}\varphi(t,\xi)|^2d\xi&=\int_{\bR^3}|M_{\delta}(t,\xi)\varphi(t,\xi)|^2d\xi\\& \lesssim \int_{\bR^3} |M_{\delta}(0,\xi )\varphi(0,\xi)|^2d\xi\\ &\lesssim\int_{\bR^3}|\langle\xi\rangle^{-4}\varphi(0,\xi)|^2d\xi.
\end{aligned}
$$
Letting $\delta \rightarrow 0$, we can obtain Theorem \ref{thm1.2} since we can take an arbitrarily large $N$.
Once the local smoothing effect of the solution is obtained, the global smoothing effect
can be reached by the same argument as in \cite{morimoto2015new}.
\end{proof}	

{\bf{Acknowledgements:}} The second author gratefully acknowledge the financial support of the National Natural Science Foundation of China (Grant No. 12301276).

\bibliographystyle{plain}

\bibliography{references}

\end{document}